\documentclass[11pt]{article}
\usepackage{amsmath,amssymb,latexsym,a4,proof}

\newtheorem{theorem}{Theorem}
\newtheorem{lemma}{Lemma}[section]

\newtheorem{corollary}[lemma]{Corollary}
\newtheorem{proposition}[lemma]{Proposition}

\newtheorem{definition}[lemma]{Definition}
\newtheorem{example}[lemma]{Example}
\newtheorem{remark}[lemma]{Remark}

\newcommand{\bl}{\begin{lemma}}
\newcommand{\el}{\end{lemma}}
\newcommand{\bt}{\begin{theorem}}
\newcommand{\et}{\end{theorem}}
\newcommand{\bcor}{\begin{corollary}}
\newcommand{\ecor}{\end{corollary}}
\newcommand{\bp}{\proof{.}}
\newcommand{\ep}{\eop}
\newcommand{\bpr}{\begin{proposition}}
\newcommand{\epr}{\end{proposition}}
\newcommand{\brem}{\begin{remark} \em}
\newcommand{\erem}{\end{remark}}
\newcommand{\bd}{\begin{definition} \em}
\newcommand{\ed}{\end{definition}}
\newcommand{\bex}{\begin{example} \em
}
\newcommand{\eex}{\end{example}}
\newcommand{\beq}{\begin{equation} }
\newcommand{\eeq}{\end{equation}}

\newcommand{\bi}{\begin{itemize}
  }
\newcommand{\ei}{\end{itemize}}
\newcommand{\ben}{\begin{enumerate} }
\newcommand{\een}{\end{enumerate} }

\newenvironment{enumr}{

\begin{enumerate}     }{\end{enumerate}

}

\newcommand{\benr}{\begin{enumr}
  }
\newcommand{\eenr}{
\end{enumr}}

\newlength{\hilflh}

\newcommand{\cL}{{\mathcal L}}

\newcommand{\cP}{{\mathcal P}}

\newcommand{\cM}{{\mathcal M}}

\newcommand{\bK}{\textbf{K}}

\newcommand{\gd}{\delta}
\renewcommand{\ge}{\varepsilon}

\newcommand{\gy}{\gamma}
\newcommand{\gw}{\omega}
\newcommand{\gS}{\Sigma}

\renewcommand{\phi}{\varphi}

\newcommand{\GL}{\mathbf{GL}}
\newcommand{\GLP}{\mathbf{GLP}}

\newcommand{\nmodels}{\nvDash}

\newcommand{\Imp}{\Rightarrow}

\newcommand{\eop}{$\Box$ \protect\par \addvspace{\topsep}}
\newcommand{\proof}[1]{\protect\par\addvspace{\topsep}\noindent {\bf Proof#1}}

\newcommand{\Kp}{\textbf{K}^+}
\newcommand{\Kxp}{\textbf{K4}^+}
\newcommand{\Sxp}{\textbf{S4}^+}
\newcommand{\Syp}{\textbf{S5}^+}
\newcommand{\tto}{\twoheadrightarrow}

\begin{document}

\title{A note on strictly positive logics and word rewriting systems}

\author{Lev Beklemishev\thanks{Steklov Mathematical Institute,
RAS; Moscow M.V. Lomonosov State University; National Research
University Higher School of Economics; email:
\texttt{bekl@mi.ras.ru}. Supported by Russian Foundation for Basic Research project 15--01--09218a and by Presidential council for support of leading scientific schools.}}

\maketitle

\begin{abstract}
We establish a natural translation from word rewriting systems to
strictly positive polymodal logics. Thereby, the latter can be
considered as a generalization of the former. As a corollary we
obtain examples of undecidable strictly positive normal modal
logics. The translation has its counterpart on the level of proofs: we formulate a natural deep inference proof system for strictly positive logics generalizing derivations in word rewriting systems. We also formulate some open questions related to the theory of modal companions of superintuitionistic logics that was initiated by L.L. Maksimova and V.V. Rybakov.
\end{abstract}

In this note we study the fragment of polymodal logic consisting of
implications of the form $A\to B$, where $A$ and $B$ are formulas
built-up from $\top$ and propositional variables using just $\land$
and the diamond modalities. We call such formulas $A$ and $B$
\emph{strictly positive} and will often omit the word `strictly.'

The interest towards such weak logics independently emerged within two different disciplines: provability logic and description logic (see~\cite{KWZ,Das12,Bek12a}). In both cases, it was observed that the strictly positive language combines simplicity and efficiency while retaining a substantial amount of expressive power of modal logic. Thus, strictly positive fragments of many standard modal logics are polytime decidable. The positive fragment of the (Kripke incomplete) polymodal provability logic $\GLP$ is both polytime decidable and complete w.r.t.\ a natural class of finite Kripke frames \cite{Das12}. The positive variable-free fragment of this logic gives rise to a natural ordinal notation system up to the ordinal $\ge_0$ and allows for a proof-theoretic analysis of Peano arithmetic~\cite{Bek12a}.

In the present paper we study some general questions related to strictly positive logics. In particular, we establish a link between proof systems for strictly positive logics and the standard word rewriting (semi-Thue) systems.

\section{Strictly positive logics}

Consider a modal language $\cL_\gS$ with propositional variables
$p,q$,\dots , a constant $\top$, conjunction $\land$, and a possibly
infinite set of symbols $\Sigma=\{a_i:i\in I\}$ understood as
diamond modalities. The family $\Sigma$ is called the
\emph{signature} of the language $\cL_\gS$. Strictly positive
formulas (or simply \emph{formulas}) are built up by the grammar:
$$A::= p \mid \top \mid (A\land B) \mid a A, \quad \text{where $a\in \Sigma$.}$$
\emph{Sequents} are expressions of the form $A\vdash B$ where $A,B$
are strictly positive formulas. We present two types of calculi for
strictly positive logics: sequent-style and deep inference-style.

Sequent-style systems for several positive logics have been
introduced and studied in \cite{Bek12a,Bek14}. This was preceded by
an equational logic characterizations of the same logics in
\cite{Das12}.

Basic sequent-style system, denoted $\Kp$, is given by the following
axioms and rules:

\ben
\item $A\vdash A; \quad A\vdash\top; \quad$ if $A\vdash B$ and $B\vdash C$ then $A\vdash C$ (syllogism);
\item $A\land B\vdash A; \quad A\land B\vdash B; \quad$ if $A\vdash B$ and $A\vdash C$ then
$A\vdash B\land C$;
\item  if $A\vdash B$ then $a A\vdash a B$.
\een

It has not been explicitly mentioned but easily follows from the
techniques of \cite{Das12,Bek14} that $\Kp$ axiomatizes the
strictly positive fragment of the polymodal version of basic modal
logic $\bK$, so we state this result witout proof.

\bt A sequent $A\vdash B$ is provable in $\Kp$ iff $\bK\vdash A \to
B$. \et

If one wishes, one can adjoin some further axioms to $\Kp$, which
correspond to some standard modal logics.

\begin{description} \item[(4)]  $aa A\vdash a A$;
\item[(T)] $A\vdash a A$;
\item[(5)] $a A\land a B\vdash a(A\land a B)$.
\end{description}

Let $\Kxp$ denote the logic axiomatized over $\Kp$ by Axiom $(4)$;
$\Sxp$ is axiomatized over $\Kp$ by $(4)$ and $(T)$; $\Syp$ is
$\Sxp$ together with $(5)$.

If $L$ is a logic, we write $A\vdash_L B$ for the statement that the
sequent $A\vdash B$ is provable in $L$. Formulas $A$ and $B$ are
called \emph{$L$-equivalent} (written $A\sim_L B$) if $A\vdash_L B$
and $B\vdash_L A$.

The following theorem is obtained by Dashkov \cite{Das12} (the case
$\Kxp$) and by Dashkov and Svyatlovsky (the cases $\Sxp$ and
$\Syp$), see~\cite{Svy14}. The latter paper also gives an infinite though explicit axiomatization
of the strictly positive fragment of the logic $\bf K4.3$.

\bt \label{frag} Let $L$ be any of the logics $\bf K4$, $\bf S4$,
$\bf S5$. Then $L\vdash A\to B$ iff $A\vdash_{L^+} B$. \et


Let $C[A/p]$ denote the result of replacing in $C$ all occurrences
of a variable $p$ by $A$. If a logic $L$ contains $\Kp$ then
$\vdash_L$ satisfies the following \emph{positive replacement
lemma}.

\bl Suppose $A\vdash_L B$, then $C[A/p]\vdash_L C[B/p]$, for any
$C$. \el \bp\ Induction on the build-up of $C$. \ep

A positive logic $L$ is called \emph{normal} if it contains $\Kp$
and is closed under the following \emph{substitution rule}: if
$A\vdash_L B$ then $A[C/p]\vdash_L B[C/p]$. It is clear that all the
positive logics considered so far are normal.

\section{Modal companions of strictly positive logics}

The language of modal logic is obtained from $\cL_\gS$ by adding
boolean connectives. Recall that a modal logic is called normal if it
contains basic modal logic $\bK$ and is closed under the rules
\emph{modus ponens}, necessitation and substitution.

There is a natural functor associating with each normal modal logic
$L$ its strictly positive fragment $\cP(L)$ consisting of all
sequents $A\vdash B$ with $A,B$ strictly positive such that $L\vdash
(A\to B)$. Vice versa, to each strictly positive normal logic $P$ we
can associate its modal counterpart $\cM(P)$ axiomatized over $\bK$
by all the implications $A\to B$ such that $A\vdash_P B$.

We note that both functors preserve inclusion, that is, are
monotone. The following obvious lemma states that $\cM$ and $\cP$,
in fact, form a \emph{Galois connection}.

\bl \label{gal} For any normal modal logic $L$ and any strictly positive normal
logic $P$,
$$\cM(P)\subseteq L \iff P\subseteq \cP(L).$$
\el

As a standard consequence we obtain that the composite operations
$\cM\cP$ and $\cP\cM$ are monotone and idempotent on the
corresponding classes of logics. Moreover, \benr
\item $\cM(\cP(L))\subseteq L$;
\item $P\subseteq \cP(\cM(P))$.
\eenr

The converse inclusions in (i) and (ii), generally, do not hold. For
(i) we can refer to the results of Dashkov~\cite{Das12}. He has
shown that for the standard modal logic $\GL$ of G\"odel and L\"ob
we have $\cP(\GL)=\Kxp$. However, by Theorem \ref{frag},
$\cM(\Kxp)={\bf K4}\neq \GL$.

For (ii), let $\Sigma=\{\Diamond\}$ and consider the logic $P$
obtained from $\Kp$ by adding the schema $\Diamond A\vdash A$. We claim that
$$\cM(P)\vdash p\land\Diamond\top\to
\Diamond p.$$ Indeed, substituting $\neg p$ for $A$ we obtain
$\cM(P)\vdash p\to \Box p$. Furthermore, $\bK\vdash \Diamond\top\land\Box p \to \Diamond p$,
therefore $\cM(P)\vdash p\land\Diamond \top\to \Diamond p$.

On the other hand, $p\land\Diamond \top\nvdash_P\Diamond p$. Consider a Kripke model
$(W,R,v)$ where $W=\{0,1\}$ and the only $R$-related elements are
$0R1$. We also let $v(p)=\{0\}$, and all the other variables are assumed to be false. For  every positive formula $A$, the set of all the nodes of $W$ where $A$ is true is downward closed. Hence, it is easy to see that this model
is sound for $P$. However $W,0\nmodels p\land\Diamond \top\to \Diamond p$.

\bigskip
It has to be noted that strictly positive logics not representable as strictly positive fragments of modal logics naturally occur in the study of reflection principles in arithmetic. For example, the system $\mathbf{RC\gw}$ axiomatizing the properties of uniform reflection principles over Peano arithmetic is of this kind~\cite{Bek14}.

A modal logic $L$ such that $\cP(L)=P$ is called a \emph{modal companion} of a positive logic $P$. As we have seen, not every normal positive logic $P$ has a companion. If it does, then $\cM(P)$ is the least modal companion of $P$ in the sense that $\cM(P)$ is contained in any other companion of $P$.
The set of modal companions of $P$, if it is not empty, also has maximal elements. This statement immediately follows from Zorn's lemma noting that the union of a chain of modal companions of $P$ is also its modal companion.

The notion of modal companion of a strictly positive logic is parallel to the one of superintuitionistic logic.  The systematic study of maximal and minimal modal companions of superintuitionistic logics was initiated by Maksimova and Rybakov \cite{MaxRyb} and followed by several important results including the Blok--Esakia theorem (see~\cite{Blo76,Esa76} and also \cite{WolZakh14} for a recent survey). For normal strictly positive logics many natural questions regarding modal companions present themselves, however so far this interesting area has not been really explored. We mention some such questions here, all of which have well-known answers in the case of superintuitionistic logics.

\begin{description}
\item[Problem 1.] Find useful criteria for a normal strictly positive logic $P$ to have a modal companion. Equivalently, for which strictly positive logics $P$ do we have $\cP(\cM(P))=P$?

\item[Problem 2.] Are there normal strictly positive logics $P$, for which there is no greatest modal companion? Are $\Kp$ and $\Kxp$ such logics?

\item[Problem 3.] Is $\GL$ a maximal modal companion of $\Kxp$? In fact, except for the cases where maximal and minimal modal companions coincide, we do not know any specific examples of maximal modal companions.
\end{description}

Let us also note that modal logics $L$ representable as the least modal companions of strictly positive logics are exactly those axiomatized over $\bK$ by a set of strictly positive implications.
Hence, if $L=\cM(P)$, as a consequence of Lemma \ref{gal} we have
$$L=\cM(P)=\cM\cP\cM(P)=\cM\cP(L).$$
Strictly positive implications are Sahlqvist formulas, therefore such logics enjoy the nice properties ensured by Sahlqvist theorem, that is, their completeness with respect to an elementary class of frames and canonicity.

Hence, we obtain the following theorem.

\bt If $L=\cM(P)$, then both $P$ and $L$ are complete w.r.t.\ an elementary class of frames. Moreover, they both are valid in the canonical frame for $L$.
\et

Obviously, an arbitrary normal strictly positive logic $P$ need not even be Kripke complete.

\section{Strictly positive deep inference calculus}

It is natural to axiomatize the consequence relation on $\cL_\gS$ in
such a way that the derived objects are positive formulas and
$A\vdash B$ is understood as provability of $B$ from hypothesis $A$.

We postulate the following conjunction introduction and elimination
rules:
$$
\infer{A\land A}{A} \qquad \infer{A}{A\land B} \qquad
\infer{B}{A\land B}
$$
The rule for $\top$ is just
$$\infer{\top}{A}$$

Notice that all the rules have one premiss. Rules in deep
inference calculi are applied within a context. A \emph{context} is
a strictly positive formula $C(p)$ in which a variable $p$ occurs
only once. Let
$$\infer{B}{A}$$ be a rule instance. For any context $C(p)$, we say that $C(B)$
is obtained from $C(A)$ by a \emph{rule application}. A
\emph{derivation} is a sequence of formulas in which every member,
except for the first one, is obtained from the previous one by a
rule application.

Let $L$ be a normal positive logic given by a set $S$ of sequents (schemata) over $\Kp$. We can naturally associate with $L$ its deep inference version $L_D$ where, in addition to the above mentioned rules for $\land$ and $\top$, for every axiom-sequent $A\vdash B$ from $S$ a rule $\displaystyle\frac{A}{B}$ in $L_D$ is postulated.

We note the following property of $L_D$.

\bl \label{contx} If $A\vdash_{L_D} B$ then $C(A)\vdash_{L_D}C(B)$,
for any context $C$. \el

\bp\ Obvious induction on the length of the derivation
$A\vdash_{L_D} B$ using the fact that if $C_1(p)$, $C_2(p)$ are
contexts then so is $C_1(C_2(p))$. \ep

\bt $A\vdash_L B$ iff $B$ is provable from $A$ in $L_D$. \et

\bp\ Both implications are established by induction on the number of rule applications in the corresponding derivation.

$(\Leftarrow)$ Assume $B$ is provable from $A$ in $L_D$ and consider the last rule application $\displaystyle\frac{B}{B'}$ in this derivation. By the induction hypothesis $A\vdash_L B'$. There is a context $C(p)$ such that $B=C(B'')$, $B'=C(A'')$ and $\displaystyle\frac{B''}{A''}$ is an instance of a postulated rule of $L_D$. Thus, we obtain  $A''\vdash_L B''$ and by positive replacement $B'=C(A'')\vdash_L C(B'')=B$. By the transitivity rule $A\vdash_L B$.

For the $(\Imp)$ we note that the syllogism rule corresponds to the
composition of derivations. The conjunction elimination axioms match
the corresponding rules.

To treat the conjunction introduction rule assume $C\vdash_L A$ and
$C\vdash_L B$. By the IH we have $L_D$ derivations of $A$ from $C$
and of $B$ from $C$. Lemma \ref{contx} yields derivations of $A\land
C$ from $C\land C$ and of $A\land B$ from $A\land C$. Hence, we can
derive in $L_D$: $C$, $C\land C$, \dots, $A\land C$, \dots, $A\land
B$, as required.

The modal rule is also interpreted by putting a deep inference proof
within a context.  If there is an $L_D$ proof of $B$ from $A$, then
by Lemma \ref{contx} there is a proof of $aB$ from $aA$. \ep

Notice that this yields deep inference systems for $\bK^+$, $\Kxp$, $\Sxp$
and $\Syp$.

\section{Word rewriting systems}

A \emph{word rewriting system} over an alphabet $\Sigma$ is given by
a set of rules of the form $A\mapsto B$ where $A,B$ are words in
$\Sigma$. Such systems are also known as \emph{semi--Thue} systems
(see \cite[Chapter 7]{DSW}). A \emph{rule application} is a substitution of an occurrence of $A$ in any word by $B$: $$XAY\to XBY.$$ A
\emph{derivation} in a system $R$ is a sequence of words in which
every member is obtained from the previous one by an application of
one of the rules of $R$. We write $A\twoheadrightarrow_R B$ iff
there is a derivation of $B$ from $A$ in $R$ (the subscript $R$ is
omitted if understood from the context).

It is well-known that finite word rewriting systems (over a finite
alphabet $\Sigma$) are a universal model of computation. In
particular, there is a finite system $R$ such that it is undecidable
whether a given word $B$ is derivable from a given word $A$.

To each word rewriting system $R$ over $\Sigma$ we associate a
normal strictly positive logic $L_R$ in $\cL_\Sigma$. $L_R$ is
obtained from $\Kp$ by adding the axioms $Ap\vdash Bp$, for each of
the rules $A\mapsto B$ from $R$. The words $A$ and $B$ are now
understood as sequences of modalities.

\bt \label{thue} $A\tto_R B$ iff $Ap\vdash_{L_R} Bp$. \et

\bp\ (only if) We argue by induction on the length of an
$R$-derivation $x$ of $B$ from $A$. Basis is easy. Suppose $x$ has
the form:
$$A\tto XUY\to XVY=B,$$ where $U\mapsto V$ is a rule from $R$. By IH
we have $Ap\vdash_{L_R} XUYp$. By the $L_R$ axiom we obtain
$UYp\vdash_{L_R} VYp$, and then by positive replacement
$XUYp\vdash_{L_R} XVYp$. Hence, $Ap\vdash_{L_R} Bp$.

The (if) part is based on the following two lemmas.

\newcommand{\lrd}{(L_R)_D}

\bl \label{top} Assume $\top$ does not occur in $A,B$. If $A\vdash_{L_R} B$ then there is a derivation of $B$ from $A$ in $\lrd$ in which the $\top$-rule is not applied. \el

\bp\ Induction on the number of applications of the $\top$-rule.
Consider any such application $$\infer{C(\top)}{C(A_1)}$$
The part of the derivation after $C(\top)$ may contain some occurrences of $\top$ inherited from this one. Replacing them all by $A_1$ yields a derivation of $B$ from $C(A_1)$ with the same number of the $\top$-rule applications. Then, the derivation $A,\dots,C(A_1),\dots, B$ has one less application of the $\top$-rule than the original derivation. \ep

\bl Assume $A,B\in \gS^*$ and $A p\vdash_{\lrd} B p$. Then there is a derivation of $B p$ from $A p$ in which no conjunction rule is applied. \el

\bp\
By Lemma \ref{top} we may assume that the $\top$-rule is not applied in the given derivation.
We argue by induction on the number of conjunction introduction
rule applications in the given derivation $d$. Since the $\top$-rule is not applied in $d$, every conjunction occurrence disappears as a result of conjunction elimination rule application either to itself, or to an external conjunction. Every formula containing at least one conjunction has the form $\gy(C_1\land C_2)$ where $\gy\in \gS^*$ (and the displayed conjunction is the outermost one).

In all the formulas of the derivation consider the outermost
conjunction. Notice that at least one outermost conjunction must be introduced in the derivation (e.g., such is the conjunction introduced first). We select the chronologically last introduced outermost conjunction. We notice that no conjunction is introduced outside this one before it is eliminated. Otherwise, the first such application would introduce an outermost conjunction later than the selected one.
Hence, the selected conjunction has exactly one successor in each step of the derivation until it disappears as a result of conjunction elimination applied to itself:
$$\gy C, \gy(C\land C), \dots ,\gd(C_1\land C_2), \gd C_i.$$
Notice that the $R$-rules do not apply to conjunctions, and the conjunction rules can only be applied inside the selected conjunction. Therefore, there exist
separate derivations of $\gd(q)$ from $\gy(q)$, and of each $C_j$ ($j\in\{1,2\}$) from $C$, respectively. It follows that we can replace this subderivation by
$$\gy C,\dots, \gd C,\dots, \gd C_i,$$ thus
eliminating at least one application of conjunction introduction
rule in the whole derivation. \ep

To complete the proof of Theorem~\ref{thue} we notice that a deep
inference format $L_R$-derivation of $B p$ from $A p$ in which no
$\top$-rule and conjunction rules are applied is essentially an $R$-derivation of $B$ from $A$. The only applicable rules are the $R$-rules
whose effect is exactly that of $R$-substitutions. \ep

\bcor There is a finitely axiomatizable undecidable strictly positive logic. \ecor

It has to be noted that the finitely axiomatized strictly positive logics that have naturally occurred so far all are polytime decidable (see \cite{Das12,Bek14}).

The results of the last section of this paper have a very close predecessor in the work of Valentin Shehtman and Alexander Chagrov (see~\cite{Sheh82,ChaShe95}). In particular, Chagrov and Shehtman exhibit undecidable propositional polymodal logics whose axioms are given by the implications of the form $A\to B$, where $A$ and $B$ are sequences of $\Box$-modalities. Clearly, such logics are the minimal modal companions of the positive logics we considered in this section. The authors, however, use semantical rather than syntactical arguments to establish a correspondence of their logics with the (semi-)Thue systems. In a sense, the correspondence between strictly positive logics of the considered kind and semi-Thue systems is even closer than for modal logics, for it extends to the level of derivations.

We thank anonymous referees for spotting some errors in the previous version of the paper and Valentin Shehtman for pointing out a connection with his work.

\bibliographystyle{plain}
\bibliography{ref-all2}

\begin{thebibliography}{10}

\bibitem{Bek12a}
L.D. Beklemishev.
\newblock Calibrating provability logic: from modal logic to reflection
  calculus.
\newblock In T.~Bolander, T.~Bra\"{u}ner, S.~Ghilardi, and L.~Moss, editors,
  {\em {Advances in Modal Logic, v.~9}}, pages 89--94. College Publications,
  London, 2012.

\bibitem{Bek14}
L.D. Beklemishev.
\newblock Positive provability logic for uniform reflection principles.
\newblock {\em Annals of Pure and Applied Logic}, 165(1):82--105, 2014.

\bibitem{Blo76}
W.J. Blok.
\newblock {\em Varieties of interior algebras}.
\newblock PhD thesis, University of Amsterdam, 1976.

\bibitem{ChaShe95}
A.~V. Chagrov and V.B. Shehtman.
\newblock Algorithmic aspects of propositional tense logics.
\newblock In {\em Lecture Notes in Computer Science}, volume 933, pages
  442--455. 1995.

\bibitem{Das12}
E.V. Dashkov.
\newblock On the positive fragment of the polymodal provability logic {GLP}.
\newblock {\em Matematicheskie Zametki}, 91(3):331--346, 2012.
\newblock English translation: \emph{Mathematical Notes} 91(3):318--333, 2012.

\bibitem{DSW}
M.~Davis, R.~Sigal, and E.J. Weyuker.
\newblock {\em Computability, complexity, and languages: fundamentals of
  theoretical computer science, 2nd ed.}
\newblock Academic Press, 1994.

\bibitem{Esa76}
L.L. Esakia.
\newblock On modal companions of superintuitionistic logics.
\newblock In {\em VII Soviet Symposium on Logic}, Kiev, 1976.

\bibitem{KWZ}
A.~Kurucz, F.~Wolter, and M.~Zakharyaschev.
\newblock Islands of tractability for relational constraints: towards dichotomy
  results for the description logic {EL}.
\newblock In {\em Advances in Modal Logic, Vol.~8}, page 271. College
  Publications, London, 2010.

\bibitem{MaxRyb}
L.L. Maksimova and V.V. Rybakov.
\newblock A lattice of normal modal logics.
\newblock {\em Algebra i Logika}, 13(2):188--216, 1974.

\bibitem{Sheh82}
V.B. Shehtman.
\newblock Undecidable propositional calculi.
\newblock In {\em Problems of cybernetics. Non-classical logics and their
  applications}, pages 74--116. Moscow, 1982.
\newblock In Russian.

\bibitem{Svy14}
M.~Svyatlovsky.
\newblock Positive fragments of modal logics.
\newblock Manuscript, in Russian.
  http://www.mi.ras.ru/${\sim}$bekl/Papers/work$\underline{\phantom{0}}$2.pdf,
  2014.

\bibitem{WolZakh14}
F.~Wolter and M.~Zakharyaschev.
\newblock {On the Blok--Esakia theorem}.
\newblock In {\em {Leo Esakia} on duality in modal and intuitionistic logics},
  volume~4 of {\em Outstanding Contributions to Logic}, pages 99--118. 2014.

\end{thebibliography}

\end{document}